\documentclass[preprint,11pt]{elsarticle}

\usepackage[margin=2.5cm]{geometry}
\usepackage{graphicx}
\usepackage{dcolumn}
\usepackage{bm}

\usepackage[utf8]{inputenc}
\usepackage[T1]{fontenc}
\usepackage{mathptmx}
\usepackage{etoolbox}
\usepackage{amsmath,amssymb}
\usepackage{color}
\usepackage{natbib}
\setcitestyle{sort}

\begin{document}

\begin{frontmatter}



\title{Finite-time self-similar rupture \\ in a generalized elastohydrodynamic lubrication model}



\author[inst1]{William Chang}
\ead{chan087@usc.edu}

\affiliation[inst1]{organization={University of Southern California},
            city={Los Angeles},
            postcode={90089}, 
            state={CA},
            country={USA}}

\author[inst2]{Hangjie Ji}
\ead{hangjie\_ji@ncsu.edu}
\affiliation[inst2]{organization={Department of Mathematics, North Carolina State University},
            city={Raleigh},
            postcode={27607}, 
            state={NC},
            country={USA}}

\begin{abstract}
Thin film rupture is a type of nonlinear instability that causes the solution to touch down to zero at finite time.
We investigate the finite-time rupture behavior of a generalized elastohydrodynamic lubrication model. This model features the interplay between destabilizing disjoining pressure and stabilizing elastic bending pressure and surface tension. The governing equation is a sixth-order nonlinear degenerate parabolic partial differential equation parameterized by exponents in the mobility function and the disjoining pressure, respectively. Asymptotic self-similar finite-time rupture solutions governed by a sixth-order leading-order equation are analyzed. In the weak elasticity limit, transient self-similar dynamics governed by a fourth-order similarity equation are also identified. 
\end{abstract}




\end{frontmatter}

\date{\today}

\section{\label{sec:intro}Introduction}
This paper presents a study of the development of finite-time singularities in a one-dimensional sixth-order partial differential equation for $h(x,t)$ on a finite domain, $0\le x \le L$,
\begin{equation}
    \frac{\partial h}{\partial t} = \frac{\partial}{\partial x}\left[h^n\frac{\partial}{\partial x}\left(
    B\frac{\partial^4 h}{\partial x^4} - \frac{\partial^2 h}{\partial x^2} + \frac{1}{mh^m}\right)\right],
    \label{eq:main}
\end{equation}
where the parameters $B, m, n > 0$.
This model is  motivated by the work by Carlson and Mahadevan \cite{carlson2016similarity} on  adhesive elastohydrodynamic touchdown that occurs as an elastic sheet begins to adhere to a wall.
The PDE \eqref{eq:main} fits into the framework of classical lubrication theory which has been widely studied for the dynamics of thin layers of slow viscous fluids spreading over solid surfaces \cite{myers,eggers2009wet}. 
Under the long-wave approximation, the lubrication equation for the evolution of the
thickness (or the height $h$ of the free-surface) of the fluid layer can be derived from Navier-Stokes equations in the low Reynolds number limit,
\begin{subequations}
\begin{equation}
    \frac{\partial h}{\partial t} = \frac{\partial}{\partial x}\left(\mathcal{M}(h) \frac{\partial p}{\partial x}\right), 
\end{equation}
where the mobility function
$
   \mathcal{M}(h) = h^n
$
with $n > 0$.
Here, $n = 3$ corresponds to the no-slip boundary condition at the liquid-solid interface, and more general Navier slip condition can be incorporated via $\mathcal{M}(h) = h^3 + \lambda h^2$. 
Following the work of Young and Stone \cite{young2017long}, we define the dynamic pressure $p$ to incorporate the elastohydrodynamic effects,
\begin{equation}
    p = B\frac{\partial^4 h}{\partial x^4} - \frac{\partial^2 h}{\partial x^2} + \Pi(h), \qquad B > 0,
\end{equation}
where ${\partial^4 h}/{\partial x^4}$ represents the elastic bending pressure due to long-wavelength sheet deformations, $B>0$ is a scaling parameter for the bending pressure, ${\partial^2 h}/{\partial x^2}$ represents the surface tension between the elastic sheet and liquid, and the disjoining pressure
\begin{equation}
   \Pi(h) = \frac{A}{h^m} , \qquad m > 0, \quad A = \frac{1}{m} >0
\end{equation}
characterizes the wetting property of the solid substrate, where  $A>0$ is the Hamaker constant.
\end{subequations}
For $m = 3$, $\Pi(h) = A/h^3$ corresponds to the van der Waals model \cite{hocking93} for the destabilizing intermolecular adhesion pressure \cite{carlson2016similarity}. Other elastohydrodynamic lubrication models \cite{hosoi2004peeling,young2017long} have also used the disjoining pressure $\tilde{\Pi}(h) = A(h^{-3}-\sigma^6h^{-9})$ with $\sigma > 0$, where the two terms in $\tilde{\Pi}(h)$ represent the repulsive and attractive intermolecular forces corresponding to the standard Lennard-Jones potential. 
A similar form of the disjoining pressure $\tilde{\Pi}(h) = Ah^{-3}(1-\epsilon h^{-1})$ with $\epsilon>0$ is often used in thin film models, setting a lower bound $h=O(\epsilon)>0$ for the film thickness $h$ and preventing thin film rupture from happening \cite{ji2017finite,bertozzi2001dewetting,schwartz2001dewetting}. 

Starting from positive and finite-mass initial data $h_0(x) > 0$ at time $t = 0$, the dynamics of the model \eqref{eq:main} are governed by the interaction between the higher-order elastic bending pressure, the surface tension, and the disjoining pressure.
Following the work of Young and Stone \cite{young2017long}, we consider the no-flux boundary conditions at $x = 0$ and $x=L$,
\begin{equation}
    h_x = h_{xxx} = h_{xxxxx} = 0 ,\quad \text{at}~ x = 0, ~L.
\label{eq:bc}
\end{equation}
The dynamics of \eqref{eq:main} can also be described by a monotone decreasing energy functional
\begin{equation}
    \mathcal{E} = \int_0^L  \frac{B}{2}\left(\frac{\partial^2 h}{\partial x^2}\right)^2+\frac{1}{2}\left(\frac{\partial h}{\partial x}\right)^2 + U(h) ~dx, 
    \qquad     \mbox{with  } \quad 
    \frac{d\mathcal{E}}{dt} = -\int_0^L h^n \left(\frac{\partial p}{\partial x}\right)^2 ~dx \le 0,
    \label{eq:energy}
\end{equation}
where $U(h)$ is the interaction potential that satisfies $U'(h) = \Pi(h)$.

Thin film rupture is a type of nonlinear instability that leads to finite-time singularities as the film thickness approaches zero at a point.
That is, $h \to 0$ at an isolated point, $x = x_c$, at a finite critical time $t = t_c$. It was shown in \cite{bernoff1998axisymmetric} that thin film equations can yield self-similar rupture singularities driven by van der Waals forces.  Different types of finite-time rupture dynamics have been investigated in a family of generalized lubrication equations parametrized by exponents in conservative and non-conservative loss terms, respectively \cite{ji2017finite,ji_witelski_2020}. In this work, we focus on the impact of the sixth-order bending pressure and the fourth-order surface tension terms on the rupture dynamics of the generalized elastohydrodynamic lubrication equation \eqref{eq:main}.

Finite-time singularities in thin film equations can result from growth in spatial perturbations due to strong instabilities. To perform a stability analysis of flat film solutions in \eqref{eq:main}, we perturb the spatially-uniform base state $h = \bar{h}$ by an infinitesimal Fourier mode
$h(x,t) = \bar{h} + \delta e^{\mathrm{i}k\pi x/L+\lambda t}+O(\delta^2)$, where $k$ is the wave number, $\lambda$ is the growth rate of disturbances, and the initial amplitude $\delta \ll 1$. Substituting the expansion into model \eqref{eq:main} and linearizing about $h = \bar{h}$ yields the dispersion relation
\begin{equation}
   \lambda = -\bar{h}^{n}\left(\frac{k\pi}{L}\right)^2\left[B\left(\frac{k\pi}{L}\right)^4+ \left(\frac{k\pi}{L}\right)^2-\frac{1}{\bar{h}^{m+1}}\right].
\label{eq:dispersion}
\end{equation}
This relation indicates that the uniform film $\bar{h} < h_c$ is long-wave unstable with respect to perturbations associated with any wave number $k \in \mathbb{Z}^+$, where the critical film thickness $h_c = \left[B({k\pi}/{L})^4+({k\pi}/{L})^2\right]^{-1/(m+1)}$.
Moreover, the relation \eqref{eq:dispersion} also shows that the disjoining pressure $\Pi(h) = 1/(mh^m)$ is destabilizing, and both the elastic bending pressure $B\partial^4h/\partial x^4$ and the surface tension $-\partial^2h/\partial x^2$ are stabilizing in the PDE \eqref{eq:main}. 

The structure of the paper is as follows. 
In Section \ref{sec:similarity} we analyze the asymptotic self-similar rupture solutions in \eqref{eq:main}, with a focus on the role of the bending pressure term. Numerical studies for the singularity solutions are presented in Section \ref{sec:numerics}, followed by concluding remarks in Section \ref{sec:conclusion}.

\section{Self-similar rupture solutions}
\label{sec:similarity}

The solutions of \eqref{eq:main} leading to rupture at a critical location $x=x_c$ for $t\to t_c$ can take the form of self-similar solutions.
Various self-similar rupture solutions of thin-film type equations have been previously analyzed \cite{carlson2016similarity, ji2017finite,witelski2000dynamics}. 
Specifically, the work of Carlson and Mahadevan \cite{carlson2016similarity} investigated the self-similar rupture solutions to a model that is equivalent to \eqref{eq:main} for $m=n=3$ without the fourth-order surface tension term.

We express the solutions of model \eqref{eq:main} using the following self-similar ansatz,
\begin{equation}
    h(x, t) \sim \tau^{\alpha}H(\eta), \quad \tau = t_c - t, \quad \eta = \frac{x-x_c}{\tau^{\beta}}, \quad \alpha, \beta > 0,
\label{eq:similarity}
\end{equation}
where the scaling parameter $\alpha > 0$ corresponds to finite-time touchdown, $h\to 0$, at $t = t_c$, and the scaling parameter $\beta$ describes the spatial focusing at $x_c$ as $\tau \to 0$. Moreover, the far-field solution $h$ away from the critical location $x_c$ should evolve slowly in time as the finite-time singularity is approached. That is, for any fixed point away from the critical location $x_c$, the time derivative term $h_t$ is bounded. This leads to the far-field boundary condition on the similarity solution $H(\eta)$,
\begin{equation}
    \alpha H - \beta\eta H_{\eta} = 0 \quad \mbox{as}~ |\eta|\to \infty.
\label{eq:farfield_bc}
\end{equation}

Substituting the ansatz \eqref{eq:similarity} into the PDE \eqref{eq:main} leads to the ordinary differential equation
\begin{equation}
\tau^{\alpha - 1}\left(-\alpha H + \eta \beta \frac{dH}{d\eta}\right)=  \tau^{n\alpha-2\beta}\frac{d}{d \eta} \left[H^{n}\frac{d}{d \eta}\left(\tau^{\alpha-4\beta}B\frac{d^4 H}{d \eta^4} - \tau^{\alpha-2\beta}\frac{d^2 H}{d \eta^2} + \frac{\tau^{-m\alpha}}{mH^m}\right)\right].
\label{eq:sim_main_PDE}
\end{equation}
For PDE models with exact similarity solutions, the values of the scaling parameters $\alpha$ and $\beta$ can be identified by separating out $\tau$ and reducing the PDE to an ODE for the similarity solution $H(\eta)$.
However, it is impossible to find an exact similarity solution for \eqref{eq:sim_main_PDE} due to the number of terms in the equation. Instead, we seek an asymptotically self-similar solution of the PDE determined by the leading-order dominant balance of terms for the limit $\tau \to 0$.

\subsection{Sixth-order similarity solution for $B=O(1)$ and $0 < n < (3m+3)/2$}
For $\tau \to 0$ with $B=O(1)$, there are four possible leading-order terms in \eqref{eq:sim_main_PDE}, the time derivative term
$\tau^{\alpha-1}(-\alpha H + \beta\eta H_{\eta})$, 
the elastic bending pressure term
$B\tau^{(n+1)\alpha-6\beta}(H^nH^{(5)})_{\eta}$,
the surface tension term
$ \tau^{(n+1)\alpha - 4\beta}(H^nH_{\eta\eta\eta})_{\eta}$, and the disjoining pressure term
$\tau^{(n-m)\alpha-2\beta}\left(H^n(\tfrac{1}{m}H^{-m})_{\eta}\right)_{\eta}$.
In the limit $\tau \to 0$,
we have $\tau^{(n+1)\alpha-6\beta} \gg \tau^{(n+1)\alpha - 4\beta}$. Therefore, the dominant balance for dynamic solutions is given by the system of equations
$
    \alpha - 1 = (n+1)\alpha - 6\beta = (n-m)\alpha - 2\beta,
$
yielding the scalings
\begin{equation}\label{eq:scaling_6}
    \alpha = \frac{2}{3m - 2n+3}, \quad \beta = \frac{m+1}{6m - 4n + 6},
\end{equation} 
and $H(\eta)$ satisfies the sixth-order similarity ODE
\begin{equation}
   -\alpha H + \beta \eta \frac{dH}{d\eta}=  \frac{d}{d \eta} \left[H^{n}\frac{d}{d \eta}\left(B\frac{d^4 H}{d \eta^4}  + \frac{1}{ mH^m}\right)\right].
 \label{eq:sim_ODE_6}
\end{equation}
With the scalings \eqref{eq:scaling_6}, the far-field boundary condition \eqref{eq:farfield_bc} reduces to 
$
    H - \frac{m+1}{4}\eta \frac{dH}{d\eta} = 0
$
as $|\eta| \to \infty$,
which indicates the asymptotic far-field behavior
$
    H(\eta)\sim C\eta^{4/(m+1)}
$
as $|\eta| \to \infty$.
The  similarity equation \eqref{eq:sim_ODE_6} corresponds to the sixth-order leading-order PDE
\begin{equation}
\frac{\partial h}{\partial t} = \frac{\partial}{\partial x}\left[h^n\frac{\partial}{\partial x}\left(
 B\frac{\partial^4 h}{\partial x^4} + \frac{1}{mh^m}\right)\right], \quad B, m, n > 0,
\label{eq:sim_PDE_6}
\end{equation}
which characterizes the balance between the sixth-order stabilizing  elastic bending pressure term and the  second-order destabilizing disjoining pressure term. Since the scaling parameters $\alpha,\beta > 0$ in \eqref{eq:scaling_6}, we need $3m-2n+3>0$, or $n < (3m+3)/2$ for the self-similar rupture solution ansatz \eqref{eq:similarity} to hold.

\subsection{Transient self-similar dynamics for $B \ll 1$ and $0<n\le m$}
\label{sec:transient}

In the weak elasticity limit, $B\ll 1$, for the regime when $B\ll \tau^{2\beta}$, the surface tension term $d^2H/d\eta^2$ dominates over the bending pressure term $Bd^4H/d\eta^4$ in \eqref{eq:sim_main_PDE}. Therefore, we have the equations of dominant balance between the time derivative term, the fourth-order stabilizing term, and the second-order destabilizing term,
$
    \alpha - 1 = (n+1)\alpha - 4\beta = (n-m)\alpha - 2\beta,
$
which leads to the scaling
\begin{equation}
    \alpha = \frac{1}{2m - n+2}, \quad \beta = \frac{m+1}{4m - 2n + 4},
\label{eq:scaling_4th}
\end{equation} 
and the similarity solution $H(\eta)$ satisfies the fourth-order nonlinear ODE
\begin{equation}
   -\alpha H + \beta \eta \frac{dH}{d\eta}=  \frac{d}{d \eta} \left[H^{n}\frac{d}{d \eta}\left( - \frac{d^2 H}{d \eta^2} + \frac{1}{ mH^m}\right)\right].
 \label{eq:4thsimilarity}
\end{equation}
The leading order terms involved represent the time derivative, the surface tension, and the disjoining pressure. 
In this case, the far-field boundary condition \eqref{eq:farfield_bc} becomes
$
    H - \frac{m+1}{2}\eta \frac{dH}{d\eta} = 0
$
as $|\eta| \to \infty$,
which indicates the asymptotic far-field behavior
$
    H(\eta)\sim C\eta^{2/(m+1)}
$
as $|\eta| \to \infty$.

The similarity equation \eqref{eq:4thsimilarity} corresponds to the fourth-order PDE
\begin{equation}
    \frac{\partial h}{\partial t} = \frac{\partial}{\partial x}\left[h^n\frac{\partial}{\partial x}\left(
 - \frac{\partial^2 h}{\partial x^2} + \frac{1}{mh^m}\right)\right], \quad m, n > 0,
    \label{eq:fourthPDE}
\end{equation}
which falls into a class of thin film-type equations studied by Bertozzi and Pugh \cite{bertozzi1998long} and
Chou and Kwong \cite{chou2007finite},
$
h_t + (h^{n}h_{xxx})_x + (h^{r}h_x)_x = 0,
$
where $n,r\in \mathbb{R}$ and $n > 0$.
This equation is identical to \eqref{eq:fourthPDE}  with $r = n-m-1$,
and the conditions for the global existence of its solutions and finite-time singularities 
have been established in terms of the exponents $n$ and $r$ in the competing second- and fourth-order terms.
Specifically, it was shown that ruptures in finite time can occur 
for the ranges $n>0$ and $r\le -1$. Therefore, the established rupture criterion for the fourth-order PDE \eqref{eq:fourthPDE} is that the exponents $m,n$ satisfy $0 < n \le m$. 
This range also guarantees that the scaling coefficients $\alpha, \beta > 0$ in \eqref{eq:scaling_4th}.

\section{Numerical studies}
\label{sec:numerics}

Next, we numerically solve the nonlinear PDE \eqref{eq:main} using a fully implicit second-order finite difference method with adaptive time stepping. The sixth-order PDE is expressed as a discretized, cell-centered system of six first-order differential equations for $h$, $k \equiv h_x$, $p \equiv k_{x}$, $q \equiv p_x$, $s \equiv Bq_x-p+\tfrac{1}{m}h^{-m}$ and $w \equiv h^n s_x$.
To identify the dynamic transition from the transient fourth-order self-similar solution to the later stage sixth-order self-similar rupture profiles, it is useful to track the relationship between the local feature of the PDE solution at the critical location $x=x_c$.
The form of the self-similar ansatz \eqref{eq:similarity} indicates that at $x = x_c$, we have
$h(x_c,t) = \tau^{\alpha} H(0)$, 
$h_{xx}(x_c,t) = \tau^{\alpha-2\beta}H''(0)$,
and $h_{xxxx}(x_c,t) = \tau^{\alpha-4\beta}H^{(4)}(0)$.
Therefore, we obtain the relation between the linearized curvature $h_{xx}(x_c,t)$ and the solution $h(x_c,t)$ at
$x=x_c$,
\begin{equation}
    h_{xx}(x_c,t) = C_1h(x_c,t)^{\nu}, \quad \text{ where } \nu=1-\frac{2\beta}{\alpha},
\label{eq:pred_hx}
\end{equation}
where the coefficient $C_1$ is uniquely determined by the local property of the similarity function $H(\eta)$. 
Based on the scaling coefficients \eqref{eq:scaling_4th} for the fourth-order self-similar dynamics and the coefficients \eqref{eq:scaling_6} for the sixth-order self-similar rupture solutions, we define the critical fourth-order and sixth-order curvature-magnitude exponent, $\nu=\nu_4(m)$ and $\nu=\nu_6(m)$, respectively,
\begin{equation}
    \nu_4(m) = -m, \quad \nu_6(m) = (1-m)/2.
\label{eq:critical_nu}
\end{equation}
We note that these critical exponents only depend on the disjoining pressure exponent $m$.
Similarly, we have
\begin{equation}
    h_{xxxx}(x_c,t) = C_2h(x_c,t)^{\mu}, \quad \text{ where } \mu=1-\frac{4\beta}{\alpha},
\label{eq:pred_hx4}
\end{equation}
which represents the relation between the elastic bending pressure and the film thickness at $x=x_c$.

\begin{figure}
    \centering
    \includegraphics[width=5.4cm]{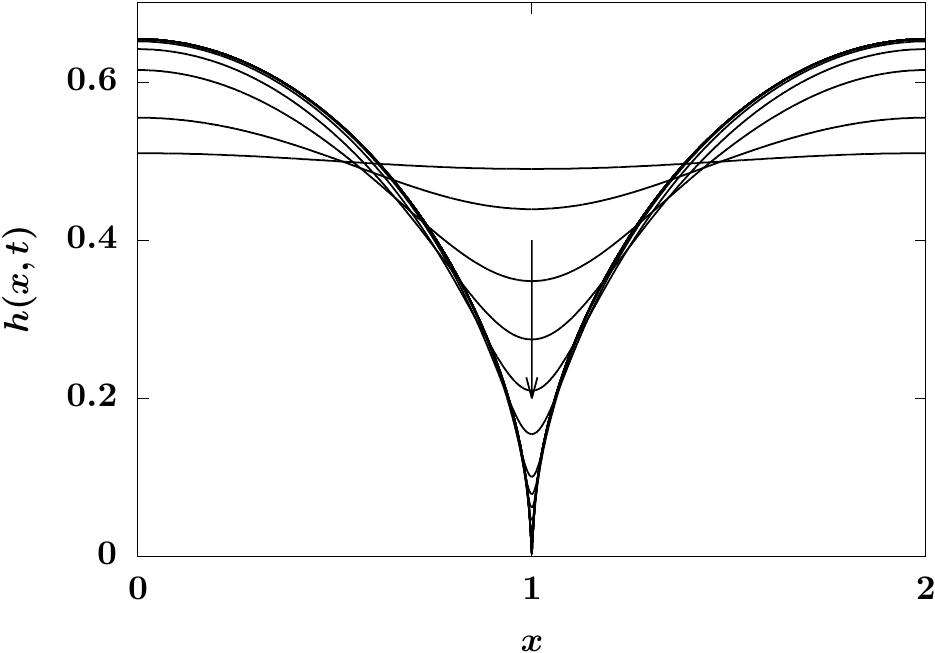}
    \includegraphics[width=5.4cm]{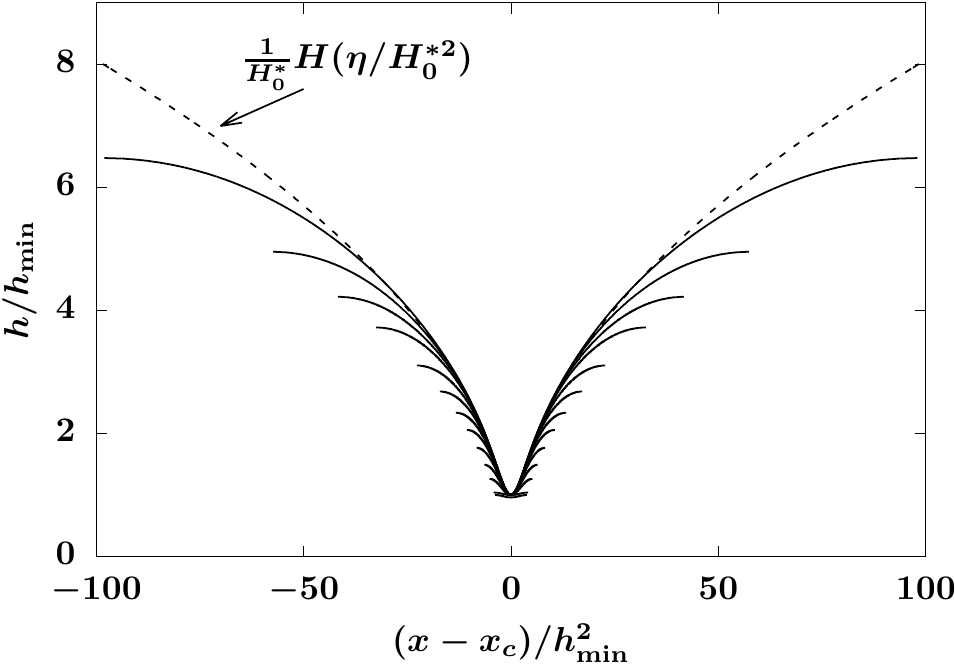}
    \includegraphics[width=5.4cm]{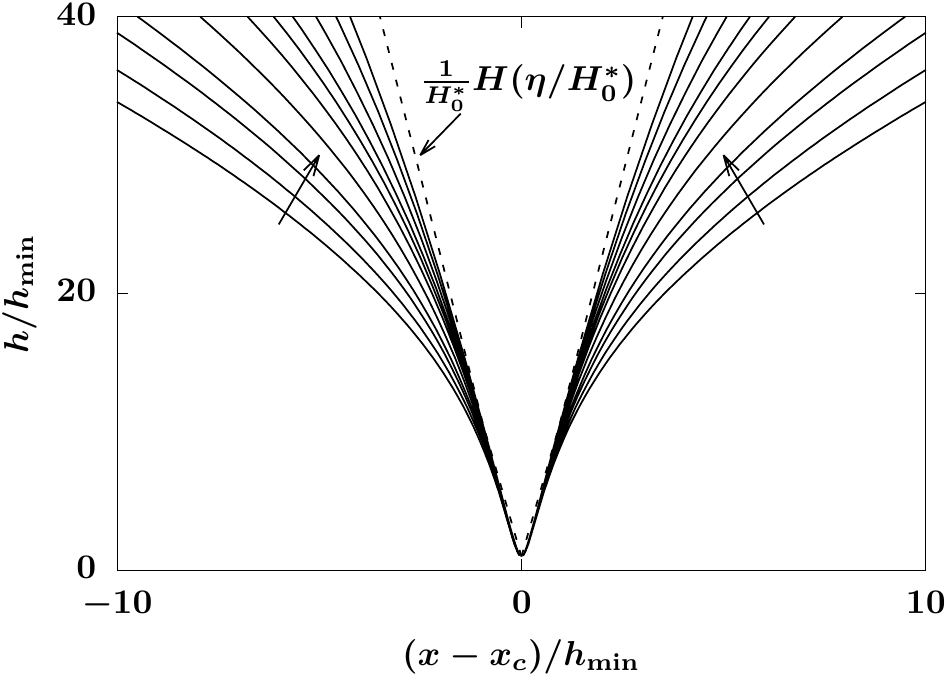} 
    \caption{(Left) Numerical solution of \eqref{eq:main} with $(m,n) = (3,3)$ and $B=10^{-5}$ and initial data $h_0(x) = 0.5+0.01\cos(\pi x)$ leading to finite-time rupture. (Center) The PDE solution for $t<0.33174$ scaled as $H(\eta/H_0^{*2})/H_0^*$ converges to the similarity solution $H(\eta)$ of the fourth-order ODE \eqref{eq:4thsimilarity}, where $H_0^*=H(0) = 0.732.$ (Right) Later stage dynamics for $0.33178<t<t_c = 0.33179$ showing that the PDE solution scaled as $H(\eta/H_0^*)/H_0^*$ converges to the similarity solution $H(\eta)$ of the sixth-order ODE \eqref{eq:sim_ODE_6}, where $H_0^*=H(0) = 2.111.$ }
    \label{fig:PDE_rup_B=1e-5}
\end{figure}

Figure \ref{fig:PDE_rup_B=1e-5} presents the dynamic solution of the PDE \eqref{eq:main} approaching a finite-time singularity at $x_c = 1$, starting from the initial condition $h_0(x) = 0.5+0.01\cos(\pi x)$ on a domain $0\le x \le 2$. This simulation corresponds to $m=n=3$ in the weak elasticity case $B = 10^{-5}$. 
For the transient self-similar stage with $B \ll \tau^{2\beta}$, the scaling parameters in the self-similar ansatz are $\alpha = 1/5$ and $\beta=2/5$ based on \eqref{eq:scaling_4th}, leading to the exponents $\nu=\nu_4 = -3$ and $\mu = -7$ in the analytical predictions \eqref{eq:pred_hx} and \eqref{eq:pred_hx4}. That is, the early-stage transient behavior satisfies $h_{xx}(x_c,t) = O(h(x_c,t)^{-3})$ and $-h_{xxxx}(x_c,t) = O(h(x_c,t)^{-7})$ at the critical location $x = x_c$. Following \cite{witelski1999stability}, we use finite difference methods to numerically solve the fourth-order similarity ODE \eqref{eq:4thsimilarity} associated with the far field boundary condition \eqref{eq:farfield_bc}
as $|\eta| \to 0$ and identify a discrete family of similarity solutions. 
Fig.~\ref{fig:PDE_rup_B=1e-5} (center) shows that the transient solution for $0<t<0.33174$, 
rescaled by $h_{\min} = \min{_x}h(x,t)$, converges to the primary similarity solution $H(\eta)$ of equation \eqref{eq:4thsimilarity} as $h_{\min}$ decreases.

\begin{figure}
    \centering
    \includegraphics[width=6cm]{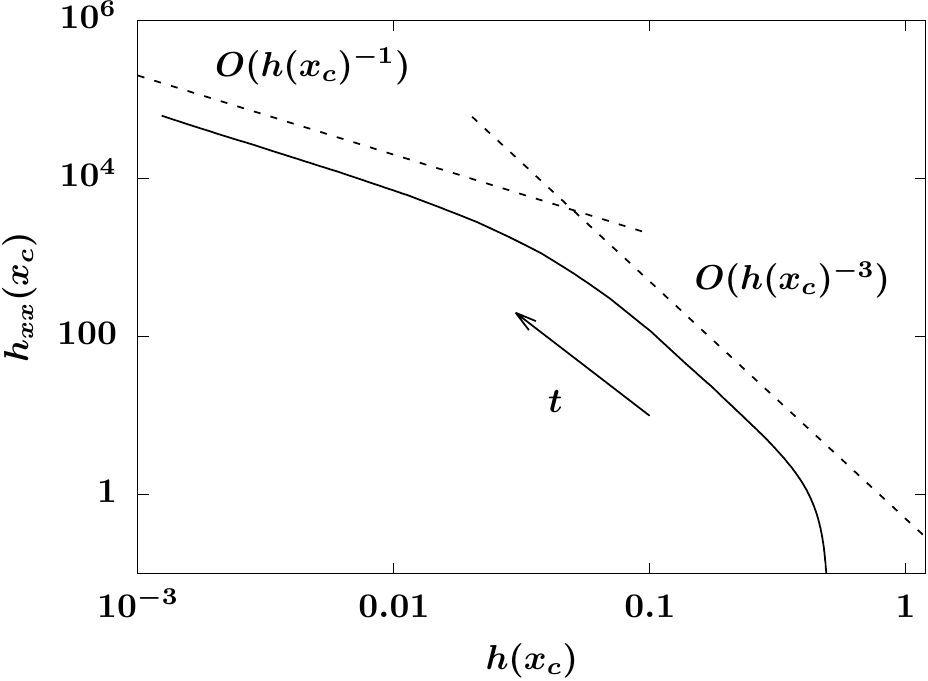}\qquad
    \includegraphics[width=6cm]{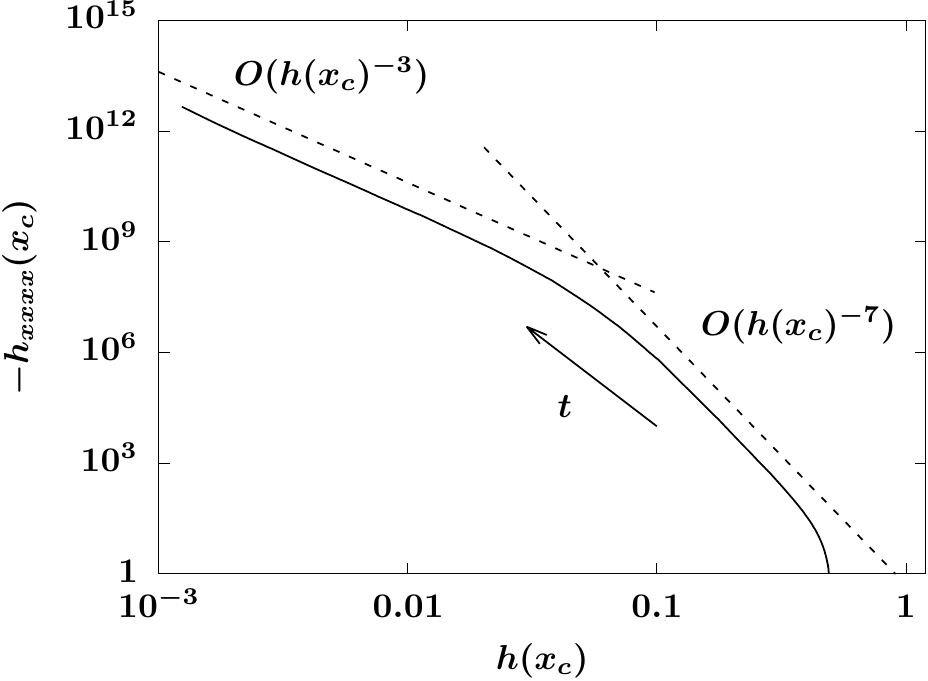}
    \caption{Plots of (left) $h(x_c,t)$ vs. $h_{xx}(x_c,t)$ and (right) $h(x_c,t)$ vs. $-h_{xxxx}(x_c,t)$ for the simulation in Fig.~\ref{fig:PDE_rup_B=1e-5}. The transition between the early stage and the later stage self-similar behaviors agree with the analytical predictions \eqref{eq:pred_hx} and \eqref{eq:pred_hx4} with the scaling parameters given by \eqref{eq:scaling_4th} for the early stage transient behavior and \eqref{eq:scaling_6} for the later stage rupture behavior.}
    \label{fig:hmin_hxx_h4_rup_B=1e-5}
\end{figure}

As the solution approaches the finite-time singularity with $\tau = t_c-t\to 0$, the condition $B\ll \tau^{2\beta}$ is no longer valid. Therefore, the PDE solution evolves following the similarity scalings \eqref{eq:scaling_6} with $\alpha = \beta = 1/3$, and the exponents in the analytical predictions \eqref{eq:pred_hx} and \eqref{eq:pred_hx4} become $\nu = \nu_6 = -1$ and $\mu = -3$, indicating that the solution satisfies $h_{xx}(x_c,t) = O(h(x_c,t)^{-1})$ and $-h_{xxxx}(x_c,t) = O(h(x_c,t)^{-3})$ as the critical time $t_c$ is approached. We plot the later stage solutions rescaled by $h_{\min}$ in Fig.~\ref{fig:PDE_rup_B=1e-5} (right) against the primary similarity solution of the sixth-order ODE \eqref{eq:sim_ODE_6}.
This transition in scaling is visible in Fig.~\ref{fig:hmin_hxx_h4_rup_B=1e-5} which depicts the relation between $h(x_c,t)$, $h_{xx}(x_c,t)$, and $h_{xxxx}(x_c,t)$ for the PDE simulation shown in Fig.~\ref{fig:PDE_rup_B=1e-5}. As $t\to t_c$, the numerically observed relations of  $h(x_c,t)$ vs. $h_{xx}(x_c,t)$ and $h(x_c,t)$ vs. $-h_{xxxx}(x_c,t)$ agree well with analytical predictions.

\begin{figure}
    \centering
    \includegraphics[width=6cm]{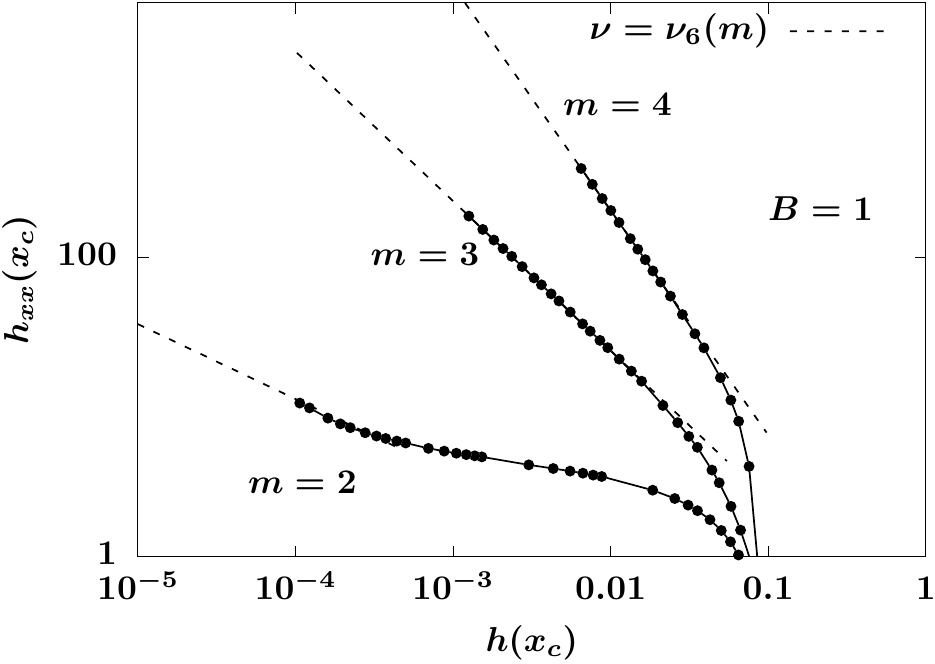}\qquad
    \includegraphics[width=6cm]{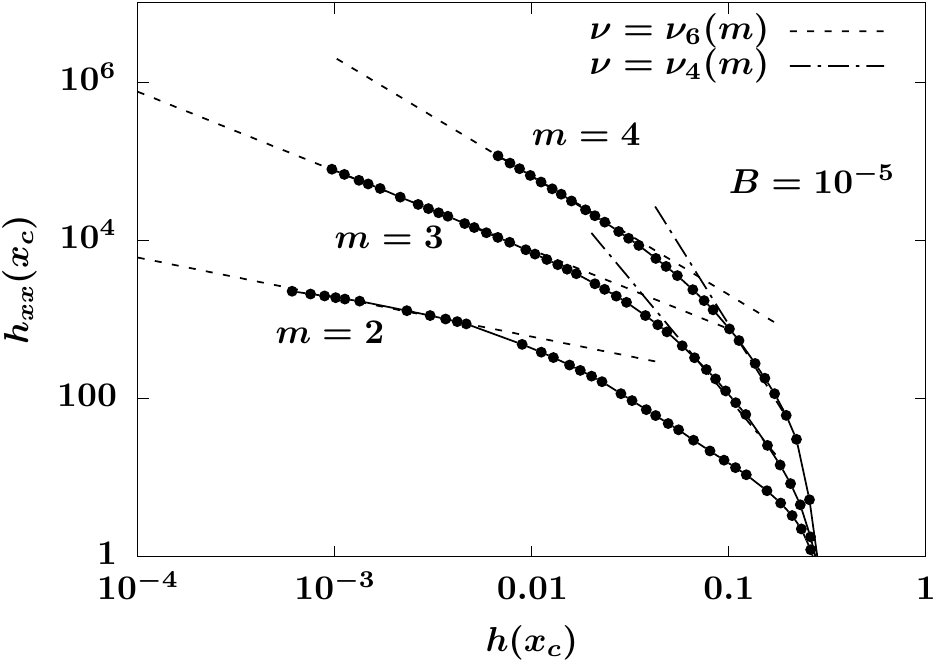}
    \caption{Dynamic solutions of PDE \eqref{eq:main}, characterized by $h(x_c,t)$ vs. $h_{xx}(x_c,t)$, starting from $h_0(x) = \bar{h}+0.01\cos(\pi x/L)$ with fixed $n = 3$, $L=0.6$ and over a range of $m$.  (Left) The strong elasticity case with $B = 1, \bar{h}=0.1$ showing sixth-order self-similar rupture behavior following the prediction \eqref{eq:pred_hx} with $\nu=\nu_6(m)$ in \eqref{eq:critical_nu}.
    (Right) The weak elasticity case with $B = 10^{-5}, \bar{h}=0.3$ showing that (a) for $m=3, 4$ that satisfies $0<n\le m$, there is a clear transition from the early-stage fourth-order self-similar dynamics with $\nu = \nu_4(m)$ to the later-stage sixth-order self-similar rupture with $\nu=\nu_6(m)$; (b) for $m=2<n$, the transient fourth-order self-similar dynamics is not observed.}
    \label{fig:bifurcation}
\end{figure}

To further investigate the rupture solution behavior and the transient dynamics in \eqref{eq:main}, 
we conduct a sequence of PDE simulations with fixed $n=3$ and over a range of $m = 2,3,4$, for both the weak ($B=10^{-5}$) and strong ($B=1$) elasticity cases. Numerical simulations starting from the initial condition $h_0(x) = \bar{h}+0.01\cos(\pi x/L)$ all lead to finite-time singularities.
To identify the rupture behaviors, we track the relation between the critical curvature $h_{xx}(x_c,t)$ and $h(x_c,t)$ in time and compare them against the predictions \eqref{eq:pred_hx} -- \eqref{eq:critical_nu}. 
Fig.~\ref{fig:bifurcation} (left) shows that in the strong elasticity case ($B=1$), the sixth-order bending pressure dominates the rupture dynamics, leading to self-similar rupture solutions following the prediction \eqref{eq:pred_hx} with $\nu = \nu_6(m)$.
Based on the discussion in Sec.~\ref{sec:transient}, 
in the weak elasticity case $B\ll 1$ with $0<n\le m$, the solution is expected to follow the transient fourth-order self-similar dynamics for $B \ll \tau^{2\beta}$ with the critical exponent $\nu_4(m)$. For the later stage dynamics towards the final rupture, the self-similar singularity occurs following the sixth-order similarity ODE \eqref{eq:sim_ODE_6} and the prediction \eqref{eq:pred_hx} with $\nu = \nu_6(m)$. 
Figure \ref{fig:bifurcation} (right) shows that in the weak elasticity case ($B=10^{-5}$), the $h(x_c,t)$ vs. $h_{xx}(x_c,t)$ curves present a clear transition from $\nu = \nu_4$ to $\nu = \nu_6$ for the cases $n\le m=3,4$. Such transition is not observed for the case  $n> m=2$, which does not satisfy the rupture criteria  \cite{chou2007finite} for the fourth-order PDE \eqref{eq:fourthPDE}. This observation confirms our analysis in Sec.~\ref{sec:transient} for the transient fourth-order rupture behavior in the weak elasticity limit.

\section{Conclusions}
\label{sec:conclusion}
This paper presents a study of the finite-time self-similar rupture dynamics in the generalized elastohydrodynamic lubrication model \eqref{eq:main} parameterized by exponents $(m,n)$ in disjoining pressure and mobility function, respectively. Asymptotically self-similar rupture solutions governed by a sixth-order nonlinear ODE are identified and numerically studied for this model. In the weak elasticity limit with $B\ll 1$ and $0<n\le m$, an interesting transition from fourth-order self-similar dynamics to the final stage sixth-order rupture solution is numerically investigated.

\section*{Acknowledgment}
H. Ji gratefully acknowledges support from Faculty Research and Professional Development Program (FRPD) provided by NC State University.


\bibliographystyle{elsarticle-num}
\biboptions{sort&compress}
{\small
\bibliography{reference}
}

\end{document}